\documentclass{article}
\usepackage{amsfonts}

\usepackage{graphicx,color,mathtools}
\usepackage{xcolor}
\usepackage{amssymb,amsmath,amsthm,mathrsfs}
\usepackage[all,cmtip]{xy}
\usepackage{epstopdf, comment, url}
\usepackage{bm} 
\usepackage{enumitem}
\usepackage{dsfont}
\usepackage{nicefrac}
\usepackage{esint}

\usepackage[pdftex,bookmarks,pdfnewwindow,plainpages=false,unicode,pdfencoding=auto]{hyperref}

\newtheorem{theorem}{Theorem}[section]
\newtheorem{corollary}[theorem]{Corollary}
\newtheorem{lemma}[theorem]{Lemma}

\theoremstyle{remark}
\newtheorem*{remark}{Remark}

\numberwithin{equation}{section}

\DeclareMathOperator{\aut}{Aut}

\DeclareMathOperator{\diam}{diam}

\DeclareMathOperator{\crit}{Crit}

\DeclareMathOperator{\hol}{Hol}

\DeclareMathOperator{\gen}{gen}

\DeclareMathOperator{\TV}{TV}

\begin{document}

\title{Invariant probability measures on holomorphic functions on the unit disk}
\author{Oleg Ivrii, Uri Kreitner}
\date{October 29, 2026}

\maketitle

\begin{abstract}
We study probability measures on the space of holomorphic functions on the unit disk that are invariant under pre-composition by conformal automorphisms and give no mass to constant functions. We determine the slowest possible growth of functions supported by such measures.

Let $M_f(R)$ denote the maximum modulus of $f$ on the hyperbolic ball of radius $R$ centered at the origin. We show that given a non-decreasing function $G:[1,\infty)\to(0,\infty)$, there exists an $\aut(\mathbb{D})$-invariant probability measure $\mu$ on $\hol(\mathbb{D}, \mathbb{C}) \setminus \{\text{constant functions} \}$ such that $\log^+ M_f(R)=O_f(G(R))$ for $\mu$ almost every $f$ if and only if $\int_1^\infty (R/G(R)) \, dR <\infty$.

For the lower bound, we use a mass transport argument to show that almost every trajectory of the gradient flow of $\log |f|$ has finite energy.  To establish sharpness, we construct holomorphic functions of slow growth on infinite Riemann surfaces.
\end{abstract}

\section{Introduction}

The complex plane $\mathbb{C}$ acts on the space of entire functions $\mathcal E$ by translations:
$$
(T_a f)(z) = f(a + z), \qquad a, z \in \mathbb{C}.
$$
Our starting point is the work of  B.~Weiss \cite{weiss}, who showed the existence of translation-invariant probability measures on $\mathcal E$ which give no mass to constant functions. A significant advance was made in \cite{BGLS}, where the authors used $\overline{\partial}$ techniques to construct invariant probability measures supported on entire functions that grow almost as slowly as possible.

In this paper, we explore analogous questions in the hyperbolic setting.
More precisely, we consider the action of $\aut(\mathbb{D})$ on the space $\hol(\mathbb{D}, \mathbb{C})$ of holomorphic functions on the unit disk by pre-composition:
$$
(T_m f)(z):= f(m(z)), \qquad m \in \aut(\mathbb{D}), \qquad  z \in \mathbb{D},
$$
and look for probability measures invariant under $\aut(\mathbb{D})$ which do not charge constant functions.
As the group $\aut(\mathbb{D})$ is not amenable, it was not clear that any such non-trivial invariant probability measures should exist.

Somewhat surprisingly, L.~Buhovsky (private communication) constructed a simple example of an
$\aut(\mathbb{D})$-invariant probability measure on $\hol(\mathbb{D}, \mathbb{C})$.
Start with a finite-area hyperbolic Riemann surface $X$ with at least one puncture. According to the uniformization theorem, one can represent $X \cong \mathbb{D}/\Gamma$ as the quotient of the unit disk by a group $\Gamma \subset \aut(\mathbb{D})$ of M\"obius transformations. By the Riemann-Roch theorem, there exists a non-constant meromorphic function $F$ on the compactification of $X$ whose poles lie at the punctures. Restricting $F$ to $X$ and lifting to the unit disk gives a holomorphic function $\widetilde{F}:\mathbb{D}\to \mathbb{C}$.

Consider the map $\widetilde{\phi}: T_1 \mathbb{D} \to \hol(\mathbb{D}, \mathbb{C})$ given by 
$$
(a, v) \mapsto \widetilde{F} \circ m_{a,v},
\qquad a \in \mathbb{D}, \quad v \in T_a \mathbb{D},
$$
where $m_{a,v} \in \aut(\mathbb{D})$ is the M\"obius transformation which takes a fixed unit tangent vector $(0, v_0)$ to $(a, v)$. 
Since $\widetilde F$ is $\Gamma$-invariant, this map descends to a map $\phi:T_1X\to\hol(\mathbb D,\mathbb C)$.
Let $\lambda_X$ denote the normalized Liouville measure on
the unit tangent bundle $T_1X$. 
By the invariance of Liouville measure, the pushforward $\phi_* \lambda_X$ is an $\aut(\mathbb{D})$-invariant probability measure on $\hol(\mathbb{D}, \mathbb{C})$ concentrated on the orbit of $\widetilde{F}$.

For a holomorphic function $f: \mathbb{D} \to \mathbb{C}$, we write
$M_f(R)$
for the maximum modulus of $f$ on a closed hyperbolic ball of radius $R$ centered at the origin.
In Buhovsky's example, $$M_{\widetilde F}(R) = \exp \bigl (\exp(R + O(1)) \bigr ).$$ In \cite{IK26}, we constructed a broad family of $\aut(\mathbb D)$-invariant probability measures supported on
Shabat functions associated to hyperbolic unimodular random planar trees of uniformly bounded degree. These functions also have at least double-exponential growth in the hyperbolic radius.
This leaves open whether invariant measures can be concentrated on functions of substantially slower growth. Our main theorem gives a sharp answer:

\begin{theorem}
\label{main-thm}
Let $G: [1, \infty) \to (0, \infty)$ be a non-decreasing function. There exists a non-trivial $\aut(\mathbb{D})$-invariant 
probability measure $\mu$ on $\hol(\mathbb{D}, \mathbb{C})$ such that for $\mu$ almost every $f$, there exist constants $C_f, R_f < \infty$ for which
$$
\log M_f(R) \le C_f \, G(R), \qquad R \ge R_f,
$$
if and only if
$$
\int_1^\infty \frac{R}{G(R)} \, dR < \infty.
$$
\end{theorem}

\begin{remark}
(i) The above criterion allows $G(R) = R^2 (1 + \log R)^{1+\varepsilon}$ for every $\varepsilon > 0$, but excludes $G(R) = R^2 (1 + \log R)$.

(ii) The mechanisms of growth in the Euclidean and hyperbolic settings are very different. In \cite{BGLS}, the growth comes from the frequent occurrence of small values and the mean value property.
In particular, $\mu$ almost every $f$ is unbounded in a unit neighbourhood of $\{|f|\le 1\}$.
In contrast, the function $\widetilde F$ in Buhovsky's example is bounded in a unit hyperbolic neighbourhood of $\{|\widetilde F|\le 1\}$.
\end{remark}

\section{A gradient flow}
\label{sec:gradient-flow}

Fix a non-constant holomorphic function $f$ on the unit disk. We consider the flow $\Phi_t = \Phi_t^f$ on 
$\mathbb{D} \setminus \crit(f)$  generated by the holomorphic vector field
\begin{equation}
\label{eq:gradient-vector-field}
V(z) = \frac{f(z)}{f'(z)} \cdot \frac{\partial}{\partial z}.
\end{equation}
Zeros of $f$ correspond to constant trajectories. Along each non-constant trajectory, the argument of $f$ remains constant while the modulus of $f$ grows exponentially:
\begin{equation}
\label{eq:gradient-flow}
f(\Phi_t(z)) = e^t f(z).
\end{equation}
In other words, the non-constant trajectories of $\Phi_t$ are the gradient lines of $\log |f|$, parametrized so that
$\frac{d}{dt}\log |f(\Phi_t(z))|=1$. 

For each starting point $z \in \mathbb{D} \setminus \crit(f)$, let
$$
\bigl (- \tau^-_f(z),\tau^+_f(z) \bigr )
$$
be the maximal interval on which the trajectory through $z$ is defined. A trajectory can have a finite lifetime only if it approaches a critical point of $f$ or leaves every compact subset of $\mathbb D$.

We call a forward trajectory $(z(t))$ {\em generic} if it is non-constant and does not converge to a critical point as $t\to\tau_f^+$. Such a trajectory leaves every compact subset of the unit disk: if $\tau_f^+ < \infty$, this follows from the continuation theorem for ODEs, while if $\tau_f^+ = \infty$, this follows from $|f(z(t))|= e^t |f(z(0))| \to\infty$. Since $f$ has countably many zeros and critical points, and only finitely many trajectories can converge to each critical point, the union of the non-generic trajectories has hyperbolic area zero.

Let $\rho(z) = \frac{2}{1-|z|^2}$ be the density of the hyperbolic metric on the unit disk. We denote the hyperbolic speed of the flow by
$$
v(z) = \rho(z) \biggl | \frac{f(z)}{f'(z)} \biggr |.
$$
For a forward trajectory $(z(t))$ starting at $z(0)$, its hyperbolic length up to time $t$ is
$$
L(t) = \int_0^t v(z(s)) \, ds.
$$
We define its energy as
$$
E_f(z(0)) = \int_0^{\tau^+_f(z(0))} v(z(s))^2 \, ds.
$$
Since $\Phi_t$ is holomorphic and carries the vector $V(z(0))$ to $V(z(t))$, its Euclidean Jacobian is
$$
\frac{|V(z(t))|^2}{|V(z(0))|^2}.
$$
Consequently, its Jacobian with respect to the hyperbolic area measure is
\begin{equation}
\label{eq:jacobian-identity}
J^f_t(z(0)) = \frac{v(z(t))^2}{v(z(0))^2}.
\end{equation}

\allowdisplaybreaks

\section{Trajectories of finite energy}

Throughout this section, $f$ will denote a non-constant holomorphic function on the unit disk. We first show:

\begin{lemma}
A generic forward trajectory of finite energy is defined for all positive time.
\end{lemma}

\begin{proof}
Suppose that the forward lifetime of a generic trajectory $(z(t))_{0 \le t < \tau^+}$ is finite. By the Cauchy-Schwarz inequality, for $0\le t<\tau^+$,
the hyperbolic distance
$$
d_h(z(0), z(t)) \, \le \, \int_0^t v(z(s))\,ds
\, \le \, \sqrt{tE_f(z(0))}
\, \le \, \sqrt{\tau^+ E_f(z(0))}.
$$
Thus $(z(t))$ is trapped in a compact hyperbolic ball, contrary to genericity.
\end{proof}

Finite energy forces the trajectory to escape slowly in the hyperbolic metric. Since the modulus of $f$ grows exponentially along the trajectory, this gives a lower bound for the growth of $f$:

\begin{lemma}
\label{finite-energy-implies-large-growth-simple}
If $f$ possesses a generic trajectory with finite energy, then
\begin{equation}
\label{eq:finite-energy-implies-large-growth-simple}
\frac{\log M_f(R)}{R^2} \to \infty.
\end{equation}
\end{lemma}

\begin{proof}
Suppose $(z(t))$ is such a trajectory. Set $a := d_h(0, z(0))$. By the previous lemma, $(z(t))$ leaves all compact subsets of the unit disk and is defined for all positive time. By the Cauchy-Schwarz inequality, for any $0 < T < t$, we have
$$
L(t) \, = \, L(T) + \int_T^t v(z(s)) \, ds \, \le \, L(T) + \sqrt{t - T} \, \biggl \{ \int_T^t v(z(s))^2 \, ds \biggr \}^{1/2},
$$
which implies that $L(t) = o(\sqrt{t})$.

For $R \ge a + 1$, let $t_R$ be the first time that the trajectory hits the circle $\partial B_h(0,R)$. Since $L(t_R) \ge R - a$,
$$
\frac{t_R}{R^2} \, \asymp \, \frac{t_R}{(R-a)^2} \, \ge \, \frac{t_R}{L(t_R)^2} \, \to \, \infty, \qquad \text{as }R \to \infty.
$$
Together with the inequality $$\log M_f(R) \, \ge \, \log |f(z(t_R))| \, = \,  \log |f(z(0))| + t_R,$$
this yields (\ref{eq:finite-energy-implies-large-growth-simple}).
\end{proof}

To obtain the integral condition in Theorem \ref{main-thm}, we will use the following consequence of Hardy's inequality:

\begin{lemma}
Suppose $v: [0, \infty) \to (0, \infty)$ is a continuous function such that
$$
E = \int_0^\infty v(t)^2 \, dt < \infty
$$
and
$$
L(t) \, = \, \int_0^t v(s) ds \, \to \, \infty, \qquad \text{as }t \to \infty.
$$
If $H: (0, \infty) \to (0, \infty)$ is a continuous function which satisfies $H(L(t)) \ge t$, then
$$
\int_0^\infty \frac{R}{H(R)} \, dR \, \le \, 2E \, < \, \infty.
$$
\end{lemma}

\begin{proof}
By Hardy's inequality \cite[Theorem 327]{HLP52} with $p = 2$,
$$
\int_0^\infty \frac{L(t)^2}{t^2} \, dt \, \le \, 4 \int_0^\infty v(t)^2 \, dt \, = \, 4E.
$$
After changing variables $R = L(t)$ and using the Cauchy-Schwarz inequality, we obtain
\begin{align*}
\int_0^\infty \frac{R}{H(R)} \, dR & = \int_0^\infty \frac{L(t)}{H(L(t))} \, v(t) dt \\
& \le \int_0^\infty \frac{L(t)}{t} \, v(t) dt \\
& \le \biggl ( \int_0^\infty \frac{L(t)^2}{t^2} \, dt \biggr )^{1/2}
\biggl ( \int_0^\infty  v(t)^2 \, dt \biggr )^{1/2} \\
& \le 2E,
\end{align*}
as desired.
\end{proof}

\begin{lemma}
\label{finite-energy-implies-large-growth}
Suppose that $f$ possesses a generic forward trajectory with finite energy. If $G: [1, \infty) \to (0, \infty)$ is a non-decreasing function such that
$$
\log M_f(R) \le C_f \, G(R), \qquad R \ge R_f,
$$
then
\begin{equation}
\label{eq:finite-energy-implies-large-growth}
\int_1^\infty \frac{R}{G(R)} \, dR < \infty.
\end{equation}
\end{lemma}

\begin{proof}
Let $(z(t))$ be such a trajectory. As in Lemma \ref{finite-energy-implies-large-growth-simple}, set $a := d_h(0, z(0))$. The function
$$
H(R) = \log M_f(R + a) - \log |f(z(0))|
$$
is positive, continuous and satisfies
$$
H(L(t)) \, \ge \, \log |f(z(t))| - \log |f(z(0))| \, = \, t.
$$
By the previous lemma,
$$
\int_1^\infty \frac{R}{H(R)} \, dR  < \infty
\qquad
\text{and hence}
\qquad
\int_{1}^\infty \frac{R}{1 + \log^+ M_f(R)} \, dR  < \infty.
$$
Together with the assumption on $G$, this shows (\ref{eq:finite-energy-implies-large-growth}).
\end{proof}

\section{Finite energy via mass transport}

Let $\mu$ be a non-trivial $\aut(\mathbb{D})$-invariant Borel probability measure on $\hol(\mathbb{D}, \mathbb{C})$. In this section, we show that for $\mu$ almost every $f$, the forward trajectory from the origin is generic and has finite energy. For this purpose, we utilize a mass transport principle of Benjamini and Schramm \cite[Theorem 5.2]{BS01}:

\begin{lemma}
\label{mass-transport-principle}
Suppose $\mathcal M$ is a positive locally finite Borel measure on $\mathbb{D} \times \mathbb{D}$ which is invariant under the diagonal action of $\aut(\mathbb{D})$. Suppose that
$$
\mathcal M(\mathbb{D} \times B) \le A_h(B), \qquad \text{for every measurable set }B \subset \mathbb{D}.
$$
Then,
$$
\mathcal M(B \times \mathbb{D}) \, = \,  \mathcal M(\mathbb{D} \times B),
$$
for every measurable set $B \subset \mathbb{D}$.
\end{lemma}

Suppose that from each point $a \in \mathbb{D}$, we send a non-negative amount of mass to every point $b \in \mathbb{D}$, with density $m(a,b) \, dA_h(a)\, dA_h(b)$. For two measurable sets $A, B \subset \mathbb{D}$, we denote by $\mathcal M(A\times B) = \int_{A\times B} m(a,b) \, dA_h(a) \, dA_h(b)$ the total mass sent from $A$ to $B$. If $m$ is invariant under the diagonal action, then all points send the same amount of mass and all points receive the same amount of mass. The above lemma then asserts that every measurable set sends exactly as much as it receives. We note that the mass transport principle also applies to singular measures.

\begin{lemma}
Let $\mu$ be a non-trivial $\aut(\mathbb{D})$-invariant Borel probability measure on $\hol(\mathbb{D}, \mathbb{C})$.
For $\mu$ almost every $f$, the trajectory $(z(t))$ starting at the origin is generic and has finite energy.
\end{lemma}

Together with Lemma \ref{finite-energy-implies-large-growth}, this shows the necessity of the integral condition in Theorem \ref{main-thm}.

\begin{proof}
Recall from Section \ref{sec:gradient-flow} that for every non-constant holomorphic function $f: \mathbb{D} \to \mathbb{C}$, the union of the non-generic trajectories has hyperbolic area zero. By the invariance of $\mu$ and Fubini's theorem, for $\mu$ almost every $f$, the trajectory starting at the origin is generic. It remains to show that the energy is finite.

For $k \in \mathbb{Z}$, let
$
S_k^f = \{ x \in \mathbb{D} : k \le \log |f(x)| < k+1 \}.
$
We will prove that
\begin{equation}
\label{eq:energy-estimate}
\mathbb{E} \biggl [ {\mathds 1}_{\{ 0 \in S_k^f \}} \frac{E_f(0)}{v_f(0)^2} \biggr ] \le 1.
\end{equation}
Since the trajectory starting at the origin is almost surely generic, $0\in S_k^f$ for some $k\in\mathbb Z$ and $0 < v_f(0) < \infty$. By (\ref{eq:energy-estimate}), $E_f(0)<\infty$ almost surely.

\medskip

{\em Step 1.}
Fix a non-constant holomorphic function $f$ and $k\in\mathbb Z$. From each $x\in S_k^f$, we send mass $J_t^f(x) dt$ to $\Phi_t^f(x)$ for $0<t<\tau_f^+(x)$. By  (\ref{eq:jacobian-identity}), the total mass sent from $x$ is
$$
\int_0^{\tau_f^+(x)} J_t^f(x) \, dt
\, = \, \frac{1}{v_f(x)^2} \int_0^{\tau_f^+(x)}v_f(\Phi_t^f(x))^2 \, dt
\, = \, \frac{E_f(x)}{v_f(x)^2}.
$$
Therefore, for measurable sets $A,B\subset\mathbb D$, the total mass sent from $A$ to $B$ is
$$
\mathcal T_f^k(A\times B)=\int_{A\cap S_k^f}\int_0^{\tau_f^+(x)}{\mathds 1}_B(\Phi_t^f(x))J_t^f(x)\,dt\,dA_h(x).
$$
We now estimate the mass received by $B$. By Fubini's theorem and the change of variables $y=\Phi_t^f(x)$,
\begin{align*}
\mathcal T_f^k(\mathbb{D} \times B)
& = \int_{S_k^f} \int_0^{\tau_f^+(x)} {\mathds 1}_{B}(\Phi_t^f(x)) J_t^f(x) \, dt \, dA_h(x) \\
& = \int_{B} \int_0^{\tau_f^-(y)} {\mathds 1}_{S^f_k}(\Phi^f_{-t}(y)) \, dt \, dA_h(y).
\end{align*}
Since $\log | f(\Phi^f_{-t}(y)) | = \log |f(y)| - t$,  the backward trajectory from $y$ spends at most one unit of time in $S_k^f$. Consequently,
$$
\mathcal T_f^k(\mathbb D\times B)\le A_h(B).
$$

{\em Step 2.} Averaging over $f$, we obtain a measure on
$\mathbb D\times\mathbb D$ given by
$$
\mathcal M_k(A\times B) = \mathbb E \mathcal T_f^k(A\times B).
$$
By Step 1, $\mathcal M_k$ is positive, locally finite and satisfies
$$
\mathcal M_k(\mathbb D\times B)\le A_h(B).
$$
Since $\mu$ is $\aut(\mathbb D)$-invariant, $\mathcal M_k$ is invariant under the diagonal
action of $\aut(\mathbb D)$ on $\mathbb D\times\mathbb D$. By the mass transport principle (Lemma \ref{mass-transport-principle}),
$$
\int_B \mathbb{E} \biggl [ {\mathds 1}_{\{x \in S_k^f\}} \frac{E_f(x)}{v_f(x)^2} \biggr ] dA_h(x) 
\, = \, \mathcal M_k(B \times \mathbb D)
\, = \, \mathcal M_k(\mathbb D \times B)
\, \le \, A_h(B).
$$
By the invariance of $\mu$, the expectation
$$
\mathbb{E} \biggl [ {\mathds 1}_{\{x \in S_k^f\}}  \frac{E_f(x)}{v_f(x)^2} \biggr ]
$$
does not depend on $x$. Taking $B$ of positive and finite area and dividing by $A_h(B)$, we obtain
$$\mathbb E \biggl [{\mathds 1}_{\{0\in S_k^f\}}\frac{E_f(0)}{v_f(0)^2} \biggr ] \le 1,$$
as desired.
\end{proof}

\section{Regularizing growth functions}

We now turn to the second half of the paper, where we construct $\aut(\mathbb{D})$-invariant measures supported on functions of slow growth. 

We say that a non-decreasing function $G: [1, \infty) \to (0,\infty)$ is a {\em regular growth envelope} if it has the following properties:
\begin{enumerate}
\item $G$ is continuous.
\item $G(1) \le 2$.
\item There exists a constant $C > 0$ such that $G(2 R) \le C \cdot G(R)$.
\item $G(R)/R^{3/2}$ is non-decreasing.
\item 
\begin{equation}
\label{eq:reciprocal-sum-z}
\int_1^\infty \frac{R}{G(R)}\, dR <\infty.
\end{equation}
\end{enumerate}

The following lemma shows that it is enough to prove the sufficiency statement in Theorem \ref{main-thm} for regular growth envelopes:

\begin{lemma}
\label{regularization-lemma}
Let $G: [1, \infty) \to(0,\infty)$ be a non-decreasing function satisfying
(\ref{eq:reciprocal-sum-z}). There exists a regular growth envelope $\widetilde G: [1, \infty) \to (0,\infty)$ with $\widetilde G\leq G$.
\end{lemma}
\begin{proof}
By the monotonicity of $G$ and \eqref{eq:reciprocal-sum-z},
$$
\sum_{k=0}^\infty \frac{4^k}{G(2^k)}<\infty.
$$
Define
\begin{equation}
\label{eq:regularization-def}
\widetilde G(R) :=  \inf_{k\ge0}
G(2^k)\max \biggl \{
\biggl(\frac{R}{2^{k+1}}\biggr)^{3/2},
\biggl(\frac{R}{2^{k+1}}\biggr)^3
\biggr \}.
\end{equation}
Choosing $k$ so that $2^k\le R<2^{k+1}$ gives
$\widetilde G(R)\le G(2^k)\le G(R)$.

We estimate the reciprocal of $\widetilde G(R)$ crudely by
\begin{equation}
\label{eq:reciprocal-bound}
\frac{1}{\widetilde G(R)}
\le
\sum_{k=0}^\infty \frac{1}{G(2^k)}
\min \biggl \{
\biggl(\frac{2^{k+1}}{R}\biggr)^{3/2},
\biggl(\frac{2^{k+1}}{R}\biggr)^3
\biggr \}.
\end{equation}
As the series on the right converges, $\widetilde G(R)>0$.

Since each term in the infimum (\ref{eq:regularization-def}) scales in the same way,
$$
a^{3/2}\widetilde G(R) \, \le \, \widetilde G(aR) \, \le \, a^3\widetilde G(R), \qquad a\ge1.
$$
Hence, $\widetilde G$ is continuous, $\widetilde G(R)/R^{3/2}$ is non-decreasing and $\widetilde G(2R)\le 8 \, \widetilde G(R)$.

A computation shows that
$$
\int_0^\infty
R\min\biggl \{
\biggl (\frac{b}{R} \biggr )^{3/2},
\biggl (\frac{b}{R} \biggr )^3
\biggr \} \, dR
 = 3b^2.
$$
Multiplying (\ref{eq:reciprocal-bound}) by $R$ and integrating, we obtain
$$
\int_1^\infty \frac{R}{\widetilde G(R)}\,dR
\, \le \,
12\sum_{k=0}^\infty \frac{4^k}{G(2^k)}
\, < \,
\infty.
$$
Finally, replacing $\widetilde{G}(R)$ by $2 \, \widetilde{G}(R) / \widetilde {G}(1)$ if necessary guarantees that $\widetilde{G}(1) \le 2$. The proof is complete.
\end{proof}

We now associate holding times $L_1,L_2,\dots$ to a regular growth envelope $G$. Set $Q(R)=G(R)/R$. Since $Q$ is continuous and strictly increasing with $Q(1) \le 2$ and $Q(R) \to \infty$ as $R \to \infty$, we can define
\begin{equation}
\label{eq:holding-times-G}
Q(r_k)=2^k,
\qquad
L_k=\lceil r_k\rceil,
\qquad k\ge1.
\end{equation}
With this definition, $2r_k \ge L_k \ge r_k \ge 1$.

\begin{lemma}
\label{properties-of-holding-times}
Let $G(R)$ be a regular growth envelope. The holding times $L_k$ given by (\ref{eq:holding-times-G}) enjoy the following properties:

{\em (i)} The sequence $(L_k)$ is non-decreasing.

{\em (ii)} $L_{k+1} \lesssim L_k$.

{\em (iii)} $L_1 + \dots + L_k \lesssim L_{k+1}$.

{\em (iv)} $G(L_k) \asymp 2^k L_k$.

{\em (v)}
\begin{equation*}
\sum_{k=1}^\infty \frac{L_k}{2^{k}} < \infty.
\end{equation*}
\end{lemma}

\begin{proof}
Property (i) follows from the monotonicity of $Q$. Since $G(R)/R^{3/2}$ is non-decreasing,
$$
Q(4 r_k) \, \ge \, 2  Q(r_k) \, = \, 2^{k+1}.
$$
Thus $r_{k+1} \le 4r_k$, and taking ceilings gives (ii).

By the doubling property, there exists an integer $m \ge 1$ such that
$$
Q(2R) < 2^m Q(R), \qquad R\ge1.
$$
Taking $R = r_k$, we get
$Q(2r_k) < 2^{m+k}$ and so $r_{k+m}>2r_k$. Since the $r_k$ are increasing, it follows that
$$
\sum_{j=1}^k r_j \lesssim r_{k+1}.
$$
The comparison $L_j \asymp r_j$ shows (iii). The same estimates yield
$$
2^{m+k} \, > \, Q(2r_k) \, \ge \, Q(L_k) \, \ge \, Q(r_k) \, = \, 2^{k},
$$
which proves (iv).

Since
$$
r_k-1
=
\int_1^\infty
{\mathds 1}_{\{Q(R)<2^{k}\}} \, dR,
$$
by Fubini's theorem,
\begin{align*}
\sum_{k=1}^\infty\frac{r_k-1}{2^k}
&=
\int_1^\infty
\sum_{k:\,Q(R)<2^{k}} \frac{1}{2^k} \, dR \\
&\le
2\int_1^\infty\frac{dR}{Q(R)}
<\infty.
\end{align*}
Together with $L_k\le(r_k-1)+2$, this proves (v).
\end{proof}

\section{An infinite Riemann surface}

Fix a regular growth envelope $G$ and its associated holding times $(L_k)$. We construct an infinite Riemann surface $X$ out of three types of pieces (caps, blocks and joins), together with a holomorphic function $F: X \to \mathbb{C}$. Caps form the bottom level, blocks preserve the boundary degree, and joins merge two columns and double the degree.

\subsection{Caps}

A {\em cap} is a copy of the bordered Riemann surface
$$
C = \bigl \{ (z, w) \, : \, w^2 = z(z-1/2)(z+1/2), \, |z| \le 1 \bigr  \}.
$$
Topologically, $C$ is a torus with one boundary component. We denote the natural projection from $C$ to $\{ |z| \le 1 \}$ by $q_C(z, w) = z$. The map $q_C$ restricts to a $2:1$ covering of the unit circle on $\partial C$. We equip each cap with the holomorphic function $F = q_C$.

\subsection{Blocks}

A {\em block} is a copy of the bordered Riemann surface
$$
B = \bigl \{ (z, w) \, : \, w^2 = z(z - e^{2/3})(z - e^{4/3}), \, 1 \le |z| \le e^2 \bigr  \}.
$$
Topologically, $B$ is a torus with two boundary components. We call the component over $\{|z|=1\}$ the {\em lower boundary} and the component over $\{|z|=e^2\}$ the {\em upper boundary}\/. The natural projection $q_B(z,w)=z$ restricts to a $2:1$ covering on each boundary component.

For $T \ge 0$ and an even integer $D \ge 2$, the function
$$
F \, = \, F_{T, D} \, = \, e^T q_B^{D/2}
$$
is $D:1$ on both components of $\partial B$. Moreover, on the lower boundary, $\log |F(z)|$ = $T$, while on the upper boundary, $\log |F(z)| = T + D$. In other words, on a block, the logarithmic modulus of $F$ increases from $T$ to $T+D$.

\subsection{Joins}

A {\em join} is a copy of the bordered Riemann surface
$$
J = \bigl \{ (z, w) \, : \, w^2 =  z - e^{1/2}, \, 1 \le |z| \le e \bigr \}.
$$
Topologically, $J$ is a pair of pants, with two boundary components over the unit circle $\{ |z| = 1 \}$ and one boundary component over $\{ |z| = e \}$. The natural projection $q_J(z,w)=z$ has degree $1$ on each lower boundary component
and degree $2$ on the upper boundary component. For $T \ge 0$ and an even integer $D \ge 2$,
$$
F \, = \, F_{T,D} \, = \, e^T q_J^D
$$
has degree $D$ on each lower boundary component and degree $2D$ on the upper boundary component. Moreover, the logarithmic modulus increases from $T$ to $T+D$.

\subsection{Assembling the surface}

We construct the Riemann surface $X$ and holomorphic function $F: X \to \mathbb{C}$ level by level. Each {\em level} consists of copies of the same piece, carrying the same function. The levels are grouped into {\em generations} with generation $0$ consisting of one level and generation $j \ge 1$ consisting of $L_j+1$ levels.

Throughout the construction, we choose the functions on the pieces so that they agree on the seams and therefore patch together to form a holomorphic function on $X$. (We glue matching boundary circles by linear dilations.) More precisely, whenever a block or a join is attached above a boundary on which $F$ has logarithmic modulus $T$ and degree $D$, we equip it with $F_{T,D}$. 

\medskip

{\em Generation $0$.} We start with countably many copies of the cap $C$, indexed by the positive integers, each carrying $F=q_C$. On their boundaries, the logarithmic modulus is $0$ and the degree is $2$.

\medskip

{\em Generation $j\ge1$.} Above each upper boundary from generation $j-1$, we place a column of $L_j$ blocks. We then pair consecutive columns and connect each pair by a join. The boundary degree remains $2^j$ throughout each column and doubles to $2^{j+1}$ at the upper boundary of the join. Consequently, across every piece of generation $j$, the logarithmic modulus increases by $2^j$. 

\medskip

We number the levels starting at $0$, with the caps forming the 0-th level. For $x\in X$, let $h(x)$ be the level of the piece containing $x$. (We assign points on a seam the {\em lower} of the two adjacent levels.)

\begin{lemma}
\label{height-and-modulus}
If $x\in X$ belongs to a piece of generation $k \ge 2$, then
\begin{equation}
\label{eq:height-and-modulus}
\log|F(x)| \asymp G(h(x)).
\end{equation}
On caps and pieces of generation 1, the function $F$ is uniformly bounded.
\end{lemma}

\begin{proof}
From the construction of $(X, F)$ and Lemma \ref{properties-of-holding-times}, it follows that the logarithmic modulus of $F$ is between
$$
\sum_{j=1}^{k-1} (L_j + 1) 2^j \qquad \text{and} \qquad \sum_{j = 1}^k (L_j + 1) 2^j
$$
and is therefore comparable to $L_k 2^k$. Similarly, the level is between
$$
\sum_{j=1}^{k-1} (L_j + 1) \qquad \text{and} \qquad \sum_{j =1}^k (L_j + 1)
$$
and is therefore comparable to $L_k$. From here, (\ref{eq:height-and-modulus}) follows from the relation $G(L_k) \asymp L_k 2^k$.
\end{proof}

\subsection{Geometry of the pieces}

While each piece has one of three possible complex structures, the hyperbolic metric it inherits from $X$ depends on how it sits in the whole surface. The following lemma provides uniform bounds:

\begin{lemma}
\label{bounded-geometry}
On each model piece, use $w$ as a local coordinate in fixed small neighbourhoods of the branch points, and $z$ elsewhere. In these coordinates, the density of the hyperbolic metric of $X$ is uniformly bounded above and below by positive constants, which do not depend on the copy of the piece.
\end{lemma}

\begin{proof}

For a piece $P \subset X$, let $P^+$ be the union of $P$ and its neighbouring pieces, so that $\overline{P}$ lies in the interior of $P^+$.
By the Schwarz lemma, $\rho_X \le \rho_{P^+}$. Since there are only three types of model pieces, $P^+$ carries one of finitely many complex structures. As $\overline{P}$ lies in the interior of $P^+$, the density of the hyperbolic metric of $X$ is uniformly bounded above on $\overline{P}$ in the stated coordinates.

For the lower bound, we note that each piece $P \subset X$ contains a curve whose free homotopy class cannot be represented by a simple curve or a power of one. Since the pieces are incompressible in $X$, the hyperbolic length of such a curve is bounded below by a universal constant. Thus, the density of the hyperbolic metric is bounded below by a universal constant at some point on each piece. By distortion estimates, this bound propagates throughout the interior of $P^+$. Since $P^+$ carries one of finitely many geometries, the bound is uniform in $\overline{P}$.
\end{proof}

\subsection[Lifting F to the universal cover]{Lifting $F$ to the universal cover}

Let $\pi: \mathbb{D} \to X$ be a universal covering map and $\widetilde F = F \circ \pi$ be a lift of $F$ to the unit disk. By Lemma \ref{bounded-geometry},
$$
|h(x_1)-h(x_2)|\le C_1 \cdot d_X(x_1,x_2) + C_2, \qquad x_1,x_2\in X.
$$
Since $\pi$ does not increase hyperbolic distances, 
$$
h(\pi(z))\le h(\pi(0)) + C_1 \cdot d_h(0,z) + C_2.
$$
Together with Lemma \ref{height-and-modulus}, this shows
\begin{equation}
\label{eq:growth-bound-ftilde}
\log^+M_{\widetilde F}(R)\lesssim G \bigl (\max \{1, h(\pi(0)) + C_1 R + C_2 \}\bigr ).
\end{equation}
Since $G$ is doubling and monotone, this simplifies to 
\begin{equation}
\label{eq:MF-bound}
\log^+ M_{\widetilde F}(R) \le C \, G(R), \qquad R \ge 1+h(\pi(0)).
\end{equation}

By Lemma \ref{bounded-geometry}, there exists an $\eta > 0$ such that every piece $P \subset X$ has hyperbolic diameter less than $\eta$.  From the construction of $F$, it is clear that the image of each piece has definite size: $\diam F(P) > c > 0$. Consequently,
\begin{equation}
\label{eq:osc-bound}
\diam \widetilde{F} (B_h(0,\eta)) > c.
\end{equation}
We emphasize that the constants $C, \eta, c > 0$ do not depend on the choice of $\pi$.

\section{Constructing invariant measures}

We now adapt Buhovsky's construction by averaging over larger and larger finite regions of $X$. For $k \ge 1$, let $X_k \subset X$ be the region that lies under a join $J_k$ of generation $k$. It contains $2^k$ caps and
$$
\underbrace{2^{k-j+1}L_j}_{\text{blocks}} \, + \, \underbrace{2^{k-j}}_{\text{joins}} \, = \, 2^k \cdot \frac{2 L_j + 1}{2^j} 
$$
pieces of generation $j$, for $1 \le j \le k$. By Lemma \ref{properties-of-holding-times}(v),
\begin{equation}
\label{eq:recurrence}
\sum_{j=1}^\infty \frac{2 L_j + 1}{2^j} < \infty.
\end{equation}
Therefore, the proportion of pieces that belong to any fixed generation converges as $k \to \infty$. Together with Lemma \ref{bounded-geometry}, the piece counts show:

\begin{corollary}
\label{tightness}
For any $\varepsilon > 0$, there exists a $j_0 \ge 1$ such that
$$
\frac{A_h(\{ a \in X_k : \gen(a) > j_0 \})}{A_h(X_k)} < \varepsilon,
$$
for all $k$. 
\end{corollary}

\begin{corollary}
\label{isoperimetry}
For every fixed $r >0$,
$$
\frac{A_h \bigl(\{ a \in X_k: d_h(a, \partial X_k)\le r \} \bigr)}{A_h(X_k)} \to 0, \qquad \text{as } k\to\infty.
$$
\end{corollary}

As in the introduction, the map  $\widetilde{\phi}: T_1 \mathbb{D} \to \hol(\mathbb{D}, \mathbb{C})$ given by 
$$
(a,v) \mapsto \widetilde{F} \circ m_{a,v} , \qquad a \in \mathbb{D}, \quad v \in T_a \mathbb{D},
$$
descends to a map $\phi: T_1 X \to \hol(\mathbb{D}, \mathbb{C})$. Let $\lambda_{X_k}$ be the restriction of the Liouville measure on $T_1 X$ to $T_1 X_k$ and set $\mu_k := (\phi_* \lambda_{X_k})/\| \phi_* \lambda_{X_k} \|_{\TV}$. We will show that $(\mu_k)$ has a subsequence converging weakly to a non-trivial $\aut(\mathbb{D})$-invariant probability measure.

For each integer $N \ge 0$, let
$$
K_N := \overline{\{\phi(a,v):(a,v)\in T_1X,\ h(a) \le N\}},
$$
where the closure is taken in the topology of locally uniform convergence.  By (\ref{eq:osc-bound}), $K_N$ contains no constant functions. By (\ref{eq:MF-bound}) and Montel's theorem, the set $K_N \subset \hol(\mathbb{D}, \mathbb{C})$ is compact and every function $f \in K_N$ satisfies 
$\log^+ M_f(R) \le C \, G(R)$ for $R \ge N+1$.
By Corollary \ref{tightness}, the sequence of measures $(\mu_k)$ is tight, and so every subsequential weak limit $\mu$ is a probability measure supported on $\bigcup_{N=0}^\infty K_N$.

Fix $m \in\aut(\mathbb D)$ and set $r := d_h(0, m(0))$. Pre-composition by $m$ takes $\widetilde F\circ m_{a,v}$ to $\widetilde F\circ m_{a,v} \circ m$, which corresponds to moving the unit tangent vector $(a,v)$ to the image of $(0,v_0)$ under $m_{a,v} \circ m$. Since this mapping moves the base point by distance at most $r$, a unit tangent vector based in $X_k$ can acquire a base point outside $X_k$ only if its original base point lies within distance $r$ of $\partial X_k$. By Corollary \ref{isoperimetry},
$$
\|\mu_k-m^*\mu_k\|_{\TV} \, \le \, 2 \cdot \frac{A_h \bigl (\{a\in X_k:d_h(a,\partial X_k) \le r\} \bigr )}{A_h(X_k)} \, \to \, 0.
$$
Since pre-composition with $m$ is continuous on $\hol(\mathbb D,\mathbb C)$, if $\mu_k \to \mu$ converges weakly, then so does $m^* \mu_k \to m^*\mu$. We conclude that every subsequential weak limit of $(\mu_k)$ is $\aut(\mathbb D)$-invariant as required.

\subsection*{Acknowledgements}

We thank Misha Sodin for posing the problem and Lev Buhovsky for sharing his insights. U.K.\ acknowledges support from the Swiss State Secretariat for Education, Research and Innovation (SERI): MB25.00004.

\subsection*{Use of AI}

The proof of Theorem \ref{main-thm} was obtained by GPT-6 Astra. The resulting argument was subsequently verified and streamlined by the authors.

\bibliographystyle{amsplain}

 \medskip

\noindent Oleg Ivrii: {\tt ivrii@tauex.tau.ac.il}

\smallskip
\noindent Wladimir Schreiber Institute of Mathematical Sciences, Tel Aviv University, Tel Aviv, Israel

\bigskip

\noindent Uri Kreitner: {\tt uri.kreitner@math.ethz.ch}

\smallskip
\noindent Department of Mathematics, ETH Z\"urich, Rämistrasse 101, 8092 Z\"urich, Switzerland
 
\end{document}